\documentclass[12pt]{article}

\usepackage{amsthm,amsmath,amsfonts}
\usepackage{latexsym}
\usepackage[T1]{fontenc}
\usepackage{a4wide}
\usepackage{times}
\usepackage{array}
\usepackage{hhline}
\usepackage{graphicx}
\numberwithin{equation}{section}

\newcommand{\ux}{\underline{x_1}}
\newcommand{\ox}{\overline{x_1}}

\newtheorem{tw}{Theorem}
\newtheorem{prop}{Proposition}
\newtheorem{lem}{Lemma}
\newtheorem{df}{Definition}
\newtheorem{fa}{Fact}

\begin{document}

\begin{center}
{\Large Steady transport equation in Sobolev-Slobodetskii spaces}
\vskip8mm
{\bf Tomasz Piasecki}
\vskip5mm
Institute of Mathematics, Polish Academy of Sciences
\vskip15mm
{\bf Abstract}
\end{center}
\vskip1mm
We show existence of a regular solution in Sobolev-Slobodetskii spaces to stationary transport equation 
with inflow boundary condition 
in a bounded domain $\Omega \subset \mathbb{R}^2$. Our result is subject to quite general 
constraint on the shape of the boundary around the points where the characteristics become 
tangent to the boundary which applies in particular to piecewise analytical domains.
Our result gives a new insight on the issue of boundary singularity for the inflow problem
for stationary transport equation, solution of which is crucial for investigation of stationary 
compressible Navier-Stokes equations with inflow/outflow.

\vskip10mm

\section{Introduction}
In this paper we are concerned with the following steady transport equation with inflow boundary condition: 
\begin{equation} \label{sys}
\begin{array}{lcr}
\sigma + U \cdot \nabla \sigma = H & \textrm{in} & \Omega,\\
\sigma = \sigma_{in} & \textrm{on} & \Gamma_{in}.
\end{array}
\end{equation}
Here, $\Omega$ is a bounded domain in $\mathbb{R}^2$, 
$\sigma$ is the unknown function and $U$ is a given vector field with normal component nonvanishing at some parts
of the boundary $\Gamma = \partial \Omega$. 
Let us denote
\begin{equation} \label{ubdr}
U \cdot n(x) = d(x),
\end{equation}
where $n$ is the outward normal to $\Gamma$. Then $\Gamma$
is divided into inflow part $\Gamma_{in}$,
outflow part $\Gamma_{out}$ and remaining part $\Gamma_0$ defined as 
\begin{equation} \label{bdr}
\begin{array}{c}
\Gamma_{in} = \{ x \in \Gamma: d(x)<0 \},\\
\Gamma_{out} = \{ x \in \Gamma: d(x)>0 \},\\
\Gamma_0 = \{ x \in \Gamma: d(x)=0 \}.
\end{array}
\end{equation}
Supposing that $H \in W^s_p(\Omega)$ and $\sigma_{in} \in W^s_p(\Gamma_{in})$,  where $W^s_p$ is fractional order Sobolev space equipped with
Slobodetskii norm (we recall the definitions below), we show existence of a solution $\sigma \in W^s_p(\Omega)$
of the system \eqref{sys} under certain assumption on $\Omega$, which turns out to be quite general and apply 
to wide class of domains.
Before we formulate the assumptions and our main result more precisely let us recall some known
results concerning steady transport equations of the type of \eqref{sys}.

One of natural and important applications of equation \eqref{sys}$_1$ is analysis of stationary compressible
Navier-Stokes system. Namely, using the concept of effective viscous flux we can eliminate divergence of
the velocity from the continuity equation which is thus reduced to steady transport equation of a type \eqref{sys}$_1$.   
This approach was applied first in \cite{BdV2} in the context of regular solutions in bounded domains
and later applied widely in the 90's, let us mention in particular 
\cite{NoPa} for an important improvement enabling to treat a problem exterior domain.
The equation \eqref{sys}$_1$ itself has also been a subject of research. 
Let us mention here results concerning regular solutions for small 
data in bounded \cite{BdV1} and exterior (\cite{No1}, \cite{No2}) domains. A stronger result on regularity in Sobolev spaces 
has been shown in \cite{GiTar}. Weak solutions for $f \in L_p, 1 \leq p \leq \infty$, $f \in \textrm{BMO}$
and $f \in {\cal H}^1$ are studied in \cite{Bae}.   
We shall underline that all above results concern the case of normal component of the velocity vanishing 
on the boundary. Admission of flow across the boundary, 
which is necessary if we want to apply above described approach to stationary compressible 
Navier-Stokes system with inflow/outflow conditions,
results in substantial mathematical difficulties
in the investigation of equation \eqref{sys}$_1$. Namely, as the equation is solved along characteristics
determined by the velocity field we have to prescribe the density on the inflow part, that is we 
obtain the boundary condition \eqref{sys}$_2$. 
The above mentioned difficulties are related to the singularity which arises
around the points where the characteristics
become tangent to the boundary. We will denote a set of these points, called briefly singularity points, by $\Gamma_s$. With the definition
\eqref{bdr} we have
\begin{equation} \label{def_gamma_s}
\Gamma_s = \overline{\Gamma_0} \cap (\overline{\Gamma_{in}} \cup \overline{\Gamma_{out}}).
\end{equation}   
The influence of boundary singularity can be seen in \cite{Kw1} where the system \eqref{sys} 
is applied in context of strong solutions to the steady compressible Navier-Stokes equations
with inflow boundary conditions.
The authors obtain a solution $u \in W^2_p$, $\rho \in W^1_p$ where $u$ is the velocity of the fluid 
and $\rho$ is the density. The result is subject to a constraint on the boundary around the singularity 
points 
which means roughly that the boundary must behave like a polynomial of order not higher then $2$,
or, in other words, that the curvature of the boundary in the singularity points must 
be strictly positive, which means that the characteristics cannot enter the domain too flatly.
There is also a limitation on the integrability of the derivatives of the solution 
$2<p<3$. It is remarkable that similar constraints on the boundary and integrability are obtained in 
\cite{PS}. These constraints motivates research towards better understanding of the singularity 
in the system \eqref{sys}. 
In \cite{TP} analogous problem with inhomogeneous slip boundary condition on the velocity
is considered in a rectangle where a problem of singularity does not appear and solutions
in the same class are obtained without limitation on $p$. The system \eqref{sys}
is solved there with a technique of elliptic regularization. In \cite{PMTP} the result
is generalized for a cylindrical domain, this time \eqref{sys} is solved with a Lagrange-type
transformation which replaces the term $u \cdot \nabla \sigma$ with a single derivative.      
Finally let us mention a recent result \cite{Ber} where the system \eqref{sys} is studied in context of weak 
solutions with inflow condition prescribing the normal component of 
$\sigma U$ instead of $\sigma$. Existence of a solution in $L_2$ is shown for right hand side in $L_2$
and a divergence free vector field $U \in H^1$. The condition on the domain is Lipschitz continuity 
without additional constraints in singularity points.  

In this paper we make a step towards better understanding of boundary singularity in system \eqref{sys}
working with strong solutions. Motivated by above described limitation in $W^1_p$ solutions
we look for solutions in slightly more general class. Namely, assuming $H \in W^s_p$ which 
is a fractional order Sobolev space with Slobodetskii norm defined below, we show existence 
of a solution $\sigma \in W^s_p$ with $s$ sufficiently small and $p$ large enough to ensure 
$sp>2$ which yields $\sigma \in L_\infty$.
Our solution is hence 'slightly worse' than a $W^1_p$ solution, 
for this price we obtain existence in a wide class of domains,
in particular piecewise analytical domains fits our assumptions.
In order to define our solutions we need a weak formulation of \eqref{sys}.
\begin{df} \label{def_sol}
By a $W^s_p$-solution to \eqref{sys} we mean a function $\rho \in W^s_p(\Omega)$ such that
\begin{equation} \label{weak}
\begin{array}{c}
\int_\Omega \sigma (\phi - U \cdot \nabla \phi - \phi {\rm div}\,U) \,dx
= -\int_{\Gamma_{in}} d \phi \sigma_{in} \,dS + \int_{\Omega} H \phi \,dx 
\quad \forall \phi \in C^1(\overline{\Omega}): \phi|_{\Gamma_{out}} = 0.
\end{array}
\end{equation}
\end{df}  
We make the following assumptions on the velocity field $U$:  
\begin{equation} \label{straight}
U = [u, 0], \quad u \in W^1_{\infty}, \quad u \geq c >0.
\end{equation}    
The assumptions \eqref{straight} are quite 
strong and definitely require a comment. 
The point is that \eqref{straight} leads to quite simple proof of the \emph{a priori}
estimate where we see easily why the approach in fractional order spaces is natural
for our problem. Therefore, introducing this assumption can be considered as a starting point 
of investigation of the problem \eqref{sys} in Sobolev-Slobodetskii spaces, which is a novelty approach. 
Extension of our result for more general fields $U$ leads to serious complications and the nature of
singularity becomes hidden by many technicalities. It will be a subject of our forthcoming paper.

{\bf Functional spaces.} We use standard Sobolev spaces $W^k_p$ with natural $k$, which consist 
of functions with the weak derivatives up to order $k$ in $L_p(\Omega)$, for the definition we refer
for example to \cite{Ad}. 
However, in view of Definition \ref{def_sol}, most important for us are Sobolev-Slobodetskii spaces 
$W^s_p$ with fractional $s$. 
For the sake of completeness we recall the
definition here. By $W^s_p(\Omega)$ we denote the space of functions for which the norm:
\begin{equation} \label{wsp}
\|f\|_{W^s_p(\Omega)} = \|f\|_{L_p(\Omega)}+\big( \iint_{\Omega^2} \frac{ |f(x) - f(y)|^p }{|x-y|^{2+sp}} \,dxdy \big)^{1/p}
\end{equation} 
is finite.
It is possible and will be convenient for us to express the $W^s_p$ - norm with a sort of fractional derivatives.
To this end let us denote the cuts of $\Omega$:
$$
\Omega_1(a) = \Omega \cap \{x_1=a\}, \quad \Omega_2(b) = \Omega \cap \{x_2=b\}.  
$$
Furthermore let us define 
\begin{equation} \label{def_ux_ox}
\ux^*={\rm inf}\{\ux(x_2): (\ux(x_2),x_2) \in \Gamma_{in}\}, \quad
\ox^*={\rm sup}\{\ox(x_2): (\ox(x_2),x_2) \in \Gamma_{out}\}.
\end{equation}
An equivalent norm in $W^s_p(\Omega)$ is defined in the following way:
\begin{fa} 
In $W^s_p$ we have an equivalent norm:
\begin{equation} \label{wsp2}
\|f\|_{W^s_p(\Omega)}^* = \|f\|_{L_p(\Omega)}+\|f\|_{W^s_p,x_1(\Omega)} + \|f\|_{W^s_p,x_2(\Omega)}, 
\end{equation}
where
\begin{equation} \label{wsp2_1}
\|f\|_{W^s_p,x_1(\Omega)} = \big( \int_a^b dx_2 \int_{\Omega_2(x_2)} dx_1 \int_{\Omega_2(x_2)} \frac{|f(y_1,x_2)-f(x_1,x_2)|^p}{|x_1-y_1|^{1+sp}} dy_1 \big)^{1/p}
\end{equation}
and
\begin{equation} \label{wsp2_2}
\|f\|_{W^s_p,x_2(\Omega)} = \big( \int_{\ux^*}^{\ox^*} dx_1 \int_{\Omega_1(x_1)} dx_2 \int_{\Omega_1(x_1)} \frac{|f(x_1,y_2)-f(x_1,x_2)|^p}{|x_2-y_2|^{1+sp}} dy_2 \big)^{1/p}.
\end{equation}
\end{fa}
\noindent
The limits of integration in outer integrals come from the definitions \eqref{def_in_out}
and \eqref{def_ux_ox}. 
The proof of equivalence of norms \eqref{wsp} and \eqref{wsp2}
can be found in \cite{Tr}.    
Let us recall also a particular case of imbedding theorem for fractional order spaces,
which will be crucial for our estimates:
\begin{fa} 
Let $f \in W^s_p(\Omega)$ with $sp>d$ where $d$ is the space dimension.
Then $f \in L_\infty(\Omega)$ and 
\begin{equation} \label{imbed}
\|f\|_{L_\infty(\Omega)} \leq C \|f\|_{W^s_p(\Omega)}.
\end{equation}    
\end{fa}   
{\bf The domain. General case.} Our result holds for a wide class of domains where the inflow and
outflow parts are defined in a natural way. A general setting is to consider inflow and outflow described as
\begin{equation} \label{def_in_out}
\begin{array}{c}
\Gamma_{in}=\{(\underline{x_1}(x_2),x_2): x_2 \in (a,b)\},\\
\Gamma_{out}=\{(\overline{x_1}(x_2),x_2): x_2 \in (a,b)\},\\ 
\end{array}
\end{equation}
with singularity points given by
$$ 
\Gamma_s = \bigcup_{1 \leq i \leq m} \{(\ux(x_2),x_2),x_2=k_{in}^i\} \; \cup \;
\bigcup_{1 \leq j \leq n} \{(\ox(x_2),x_2),x_2=k_{out}^j\} 
$$
where 
$$
a=k^1_{in} < \ldots < k^m_{in}=b, \quad a=k^1_{out} < \ldots < k^n_{out}=b. 
$$
Taking into account the definition of singularity points \eqref{def_gamma_s} 
and assumption \eqref{straight}, around the singularity points we must have
$$
{\rm lim}_{x_2 \to k_{in}^{i-}} |\ux'(x_2)| = {\rm lim}_{x_2 \to k_{in}^{i+} } |\ux'(x_2)| = \infty,
$$$$
{\rm lim}_{x_2 \to k_{out}^{i-}} |\ox'(x_2)| = {\rm lim}_{x_2 \to k_{out}^{i+}}  |\ox'(x_2)|  = \infty.
$$
%
%
Furthermore we assume 
\begin{equation} \label{def_sing3}
\begin{array}{c}
{\rm lim}_{x_2 \to a^+} \underline{x_1}(x_2) = c_1, \quad {\rm lim}_{x_2 \to b^-} \underline{x_1}(x_2) = c_2,\\
{\rm lim}_{x_2 \to a^+} \overline{x_1}(x_2) = d_1, \quad {\rm lim}_{x_2 \to b^-} \overline{x_1}(x_2) = d_2,
\end{array}
\end{equation}
with $c_i \leq d_i$. Then we have 
\begin{equation}
\Gamma_0 = (c_1,d_1) \times \{a\} \cup (c_2,d_2) \times \{b\} \cup \Gamma_s.
\end{equation}
   
{\bf The domain. Simple representative case.} 
As we shall see, the proof relies essentially on investigation of the behaviour of the boundary near the singularity points.
Provided that around each singularity point condition \eqref{flat} below is satisfied,  
it is enough to consider a following simple domain with
inflow and outflow parts of the boundary given by \eqref{def_in_out}, \eqref{def_sing3} with 
$a=c_i=d_i=0$
and
\begin{equation} \label{def_sing2}
\begin{array}{c}
{\rm lim}_{x_2 \to 0^+} \underline{x_1}'(x_2) = {\rm lim}_{x_2 \to b^-} \overline{x_1}'(x_2)= -\infty,\\
{\rm lim}_{x_2 \to 0^+} \overline{x_1}'(x_2) = {\rm lim}_{x_2 \to b^-} \underline{x_1}'(x_2)= +\infty.
\end{array}
\end{equation}
Then we have two singularity points:
$$
\Gamma_s=\Gamma_0=\{(0,0),(0,b)\}.    
$$
%
We assume further that these are the only singularity points, that is, there are no singularity points 
'inside' $\Gamma_{in}$ and $\Gamma_{out}$. 
%
Similar domain is considered in \cite{Kw1}, hence our choice facilitates a comparison of the results 
and techniques with above mentioned paper. 
Around each singularity point $x_2$
is given as a function of $x_1$. We assume that this function satisfies the following constraint:
\begin{equation} \label{flat}  
|x_2(x_1) - x_2(y_1)| \geq C |x_1-y_1|^r \quad \textrm{for some} \quad r>1. 
\end{equation}
This condition is crucial in our proof and it deserves a comment. 
Condition \eqref{flat} means that the boundary around the singularity points 
is allowed to behave in particular like a polynomial of arbitrary degree.
It is worth to compare our constraint with \cite{Kw1} where analogous condition 
is required with $r=2$, therefore our restriction is much more general. 
The condition (\ref{flat}) is quite technical, however it applies for a wide class of functions,
in particular it holds true if the boundary around the singularity points is an analytic function.
It is quite basic result but we show it for the sake of completeness.
\begin{lem}
Assume that $x_2$ is an analytic function of $x_1$ around the singularity points. 
Then (\ref{flat}) holds.
\end{lem}
\noindent   
\emph{Proof.} 
By \eqref{def_sing2} we have $x_2'(x_1)=0$ at the singularity points. Therefore 
it is enough to show that if $f: \mathbb{R} \to \mathbb{R}$ is analytic in some $[-r,r]$,
$f(0)=f'(0)=0$ and $f \neq 0$ then
\begin{equation}
|f(x)| \geq C |x|^N \quad \textrm{for} \quad x \in (-l,l)
\end{equation}   
for some $C>0$, $N \geq 2$ and $l<r$ sufficiently small.  
Since $f \neq 0$ and $f$ is analytic, we must have $f^{(n)}(0) \neq 0$ for some $n \geq 2$.
Let $f^{(k)}(0)$ be the first derivative not vanishing in $0$. Then we have
$$
f(x) = \frac{f^{(k)}(0)}{k!} x^k + R^{k+1}(x),
$$
where $|R^{k+1}(x)| \leq M |x|^{k+1}$ for $x \in (-r,r)$.
Hence 
$$
|f(x)| \geq \Big|\frac{f^k(0)}{k!}\Big| |x|^k-M|x|^{k+1} = \Big( \frac{f^k(0)}{k!} - M|x|\Big) |x|^k. \square
$$ 

{\bf Main result.} We are now ready to formulate our main result. 
\begin{tw} \label{main} 
Assume that the boundary of $\Omega$ around the singularity points satisfy
\eqref{flat} for some $C>0$ and $N \in \mathbb{N}$.
Assume $U$ satisfy the assumptions \eqref{straight}. Assume further that 
$H \in W^s_p(\Omega)$ and $\sigma_{in} \in W^s_p(\Gamma_{in})$ for $s,p$ such that
$\frac{1}{r}>s>\frac{2}{p}$ where $r$ is the exponent from \eqref{flat}. 
Then there exists a solution $\sigma$ to \eqref{sys} such that 
\begin{equation} \label{est_trans}
\|\sigma\|_{W^s_p(\Omega)} \leq C [\|H\|_{W^s_p(\Omega)} + \|\sigma_{in}\|_{W^s_p(\Gamma_{in})}].
\end{equation}
\end{tw}
   
\section{Proof of Theorem \ref{main}}
As we are concerned with linear system, the core of the proof is in appropriate estimates.
The following proposition gives \emph{a priori} estimate in $W^s_p$ for a solution of 
\eqref{sys}.
\begin{prop}  
Assume $\Omega$, $H$, $U$ and $\sigma_{in}$ satisfy the assumptions of Theorem \ref{main}. 
Let $\sigma$ be a sufficiently regular solution to the equation \eqref{sys}. Then 
\eqref{est_trans} holds. 
\end{prop}
\noindent
\emph{Proof.} First of all, let us notice that it is enough to show \eqref{est_trans} 
for the solution of the equation 
\begin{equation} \label{sys2}
\begin{array}{c}
\tilde \sigma_{x_1} = \tilde H \in W^s_p(\Omega), \nonumber \\[3pt]
\tilde \sigma|_{\Gamma_{in}} =\tilde \sigma_{in} \in W^s_p(\Gamma_{in}) .
\end{array}
\end{equation}
Indeed, if we define 
\begin{displaymath}
V(x)=\int_{\ux(x_2)}^{x_1}\frac{1}{u(s,x_2)}\,ds,  \quad  \tilde \sigma = e^V \sigma,
\end{displaymath}
then we have
\begin{equation} \label{dx1}
\partial_{x_1}\tilde \sigma = V_{x_1} e^V\sigma + e^V\sigma_{x_1}.
\end{equation}
%
Since $V_{x_1}=\frac{1}{u}$, multiplying \eqref{dx1} by $u$ we get 
$$
u\tilde \sigma_{x_1} = [\sigma + u\sigma_{x_1}]e^V
=He^V,
$$
so
\begin{equation}
\tilde \sigma_{x_1} = \frac{He^V}{u}=:\tilde H.
\end{equation}
Now due to boundedness of $\Omega$ and assumptions on $u$ we have $\frac{e^V}{u} \in W^1_{\infty}$ and
$0<M_1 \leq \frac{e^V}{u}\leq M_2$ in $\Omega$, therefore 
\begin{equation}
C_1 \|\frac{e^V}{u} f\|_{W^s_p} \leq \|f\|_{W^s_p} \leq C_2 \|\frac{e^V}{u} f\|_{W^s_p}
\end{equation} 
for $f \in W^s_p(\Omega)$, which shows equivalence of estimates for \eqref{sys2} and \eqref{sys}. $\square$  

Let us proceed with the proof of Theorem 1 for \eqref{sys2} (skipping tildes for simplicity).
For convenience we use the norm $\|\cdot\|_{W^s_p}^*$ (\ref{wsp2}), hence 
it is enough to find the bounds on the norms $\|\sigma\|_{L_p}$,
$\|\sigma\|_{W^s_p,x_1}$ and $\|\sigma\|_{W^s_p,x_2}$ defined in (\ref{wsp2_1}) and (\ref{wsp2_2}).
The idea is to express pointwise values of $\sigma$ with integrals of $H$ along the characteristics
of \eqref{sys}$_1$, that is, straight lines due to assumption \eqref{straight}.
In fact, to operate with pointwise values we need sufficiently smooth functions, therefore we
consider smooth approximations of $\sigma$ and $H$ and use standard density argument. 
We start with estimate for $\|\sigma\|_{L_p}$. We have
$$
\sigma(x_1,x_2)=\sigma_{in}(\ux(x_2),x_2)+\int_{\ux(x_2)}^{x_1} H (t,x_2)dt,
$$
therefore we directly get 
\begin{equation} \label{e0}
\|\sigma\|_{L_p(\Omega)} \leq C [\|\sigma_{in}\|_{L_p(\Gamma_{in})}+\|H\|_{L_p(\Omega)}].
\end{equation}
Now we consider $\|\sigma\|_{W^s_p,x_2(\Omega)}$. For convenience denote $h:=y_2-x_2$. 
We can assume that $\ux(y_2)<\ux(x_2)$, otherwise we interchange $\ux(y_2)$ and $\ux(x_2)$ in the integrals. 
We write
$$
\sigma(x_1,x_2+h) - \sigma(x_1,x_2) = \sigma_{in}(\ux(x_2+h))-\sigma_{in}(\ux(x_2)) + 
$$$$
\int_{\underline{x_1}(x_2+h)}^{x_1} H(t,x_2+h) \,dt - 
\int_{\underline{x_1}(x_2)}^{x_1} H(t,x_2) \,dt
$$$$
= \sigma_{in}(\ux(x_2+h))-\sigma_{in}(\ux(x_2)) +
 \int_{\underline{x_1}(x_2)}^{x_1}  [H(t,x_2+h) - H(t,x_2)] \,dt \
+ \int_{\underline{x_1}(x_2+h)}^{\underline{x_1}(x_2)} H(t,x_2+h) \,dt =
$$$$
=: I_0 + I_1 + I_2,
$$
Hence $|\sigma(x_1,x_2+h) - \sigma(x_1,x_2)|^p \sim |I_0|^p + |I_1|^p + |I_2|^p$. 
We omit the limits of integration w.r.t. $x_1$ and $x_2$ since they do not play role in the computations.
Concerning the limits of $h=y_2-x_2$, these depends on $x$, but what is important is integrability around
$0$, therefore we assume without loss of generality $h>0$ and integrate 
with respect to $h$ from $0$ to some $\delta>0$.  
First of all, notice that we have $dx_2 \leq dS(x_2)$ where $dS(x_2)$ is the boundary measure at $\Gamma_{in}$.
Therefore
\begin{align} \label{e10}
\int dx_1 \int dx_2 \int_0^{\delta} \frac{|I_0|^p}{|h|^{1+sp}} \,dh & \leq 
\int dx_1 \iint_{\Gamma_{in}^2} \frac{ |\sigma_{in}(\ux(y_2))-\sigma_{in}(\ux(x_2))|^p }{|x_2-y_2|^{1+sp}} dS(x_2)dS(y_2) \nonumber \\[5pt]
& \leq C(\Omega) \|\sigma_{in}\|_{W^s_p(\Gamma_{in})}.
\end{align} 
Next, by Jensen inequality we have
$$
\int dx_1 \int dx_2 \int_0^\delta \frac{|I_1|^p}{|h|^{1+sp}} \,dh = \int dx_1 \int dx_2 \int_0^\delta dh \frac{| \int_{\underline{x_1}(x_2)}^{x_1}  [H(t,x_2+h) - H (t,x_2)] \,dt |^p}{|h|^{1+sp}} 
$$$$
\leq \int dx_1 \int dx_2 \int_0^{\delta} dh \frac{ |x_1-\underline{x_1}(x_2)|^{p-1} \int_{\underline{x_1}(x_2)}^{x_1} |H(t,x_2+h)-H(t,x_2)|^p \,dt } {|h|^{1+sp}}
$$$$
\leq C(\Omega) \int dx_1 \int dx_2 \int_{\underline{x_1}(x_2)}^{x_1} \int_0^{\delta} \frac{|H(t,x_2+h)-H(t,x_2)|^p}{|h|^{1+sp}} \,dh\,dt
$$$$
\leq C(\Omega) \int dx_1 \int dx_2 \int_{\underline{x_1}(x_2)}^{\overline{x_1}(x_2)} \int_0^{\delta} \frac{|H(t,x_2+h)-H(t,x_2)|^p}{|h|^{1+sp}} \,dh\,dt.
$$
which yields 
\begin{equation} \label{e11}
\int dx_1 \int dx_2 \int_0^\delta \frac{|I_1|^p}{|h|^{1+sp}} \,dh \leq C(\Omega)\|H\|_{W^s_p,x_2}^p.  
\end{equation}
Note that the above estimate did not involve any assumption on the boundary. The assumption \eqref{flat} 
will play role in the second part of the estimate. Namely, for $I_2$ we have
$$
\int dx_1 \int dx_2 \int_0^{\delta} \frac{|I_2|^p}{|h|^{1+sp}}\,dh = 
\int dx_1 \int dx_2 \int_0^{\delta} \frac{ \big|\int_{\underline{x_1}(x_2+h)}^{\underline{x_1}(x_2)} H(t,x_2+h)^p \,dt\big| }{|h|^{1+sp}} \,dh
$$$$
\leq \int dx_1 \int dx_2 \int_0^{\delta} dh \frac{ |\underline{x_1}(x_2)-\underline{x_1}(x_2+h)|^{p-1} \int_{\underline{x_1}(x_2+h)}^{\underline{x_1}(x_2)}|H(t,x_2+h)|^p \,dt}{|h|^{1+sp}}    
$$$$
\leq \int dx_1 \int dx_2 \int_0^{\delta} dh  \frac{ \|H\|_{L_{\infty}}^p |\underline{x_1}(x_2)-\underline{x_1}(x_2+h)|^p }{|h|^{1+sp}}.
$$
Now condition (\ref{flat}) implies that 
\begin{equation} \label{flat_inv}
|\underline{x_1}(x_2) - \underline{x_1}(x_2+h)| \leq C |h|^\epsilon
\end{equation}
where $\epsilon=\frac{1}{r}$, therefore 
$$
\int dx_1 \int dx_2 \int_0^{\delta} \frac{|I_2|^p}{|h|^{1+sp}}\,dh \leq 
\|H\|_{L_\infty}^p \int dx_1 \int dx_2 \int_0^\delta |h|^{-1+p(\epsilon-s)}. 
$$
The last integral is finite provided that $s<\epsilon$.
By the imbedding theorem we get  
\begin{equation}
\int dx_1 \int dx_2 \int_0^{\delta} \frac{|I_2|^p}{|h|^{1+sp}}\,dh \leq C \|H\|_{W^s_p}^p,  \label{e12}
\end{equation}  
provided that $s<\epsilon$ where $\epsilon$ is the exponent from (\ref{flat_inv}).
On the other hand we require $s>\frac{2}{p}$.
Therefore we can take $p>2r$ where $r$ is from \eqref{flat} and 
$\frac{1}{r}>s>\frac{2}{p}$.
From \eqref{e10}, \eqref{e11} and \eqref{e12} we conclude 
\begin{equation} \label{e1}
\|\sigma\|_{W^s_p,x_2(\Omega)} \leq C [\|H\|_{W^s_p(\Omega)} + \|\sigma_{in}\|_{W^s_p(\Gamma_{in})}].
\end{equation}

It remains to estimate the fractional norm with respect to $x_1$. The estimate is more direct
and does not involve any assumptions on the geometry of the boundary. We have
$$
\sigma(x_1+h,x_2)-\sigma(x_1,x_2) = \int_{x_1}^{x_1+h} H(t,x_2) \,dt.
$$ 
Now
$$
\int dx_1 \int dx_2 \int_0^{\delta} dh \frac{ |\int_{x_1}^{x_1+h} H(t,x_2) dt|^p }{|h|^{1+sp}} \leq
$$$$
\leq \int dx_1 \int dx_2 \int_0^\delta dh \frac{|h|^{p-1}\int_{x_1}^{x_1+h} |H(t,x_2)|^p \,dt}{|h|^{1+sp}}
\leq C \|H\|_{L_\infty}^p \int_0^\delta |h|^{-1+p(1-s)} \,dh.
$$
The last integral is finite for $s<1$. We see that
\begin{equation} \label{e2}
\|\sigma\|_{W^s_p,x_1} \leq C \|H\|_{W^s_p}
\end{equation}  
for $H \in W^s_p$ where $1>s>\frac{2}{p}$, what is clearly weaker assumption then in the previous estimate.
Combining \eqref{e0}, \eqref{e1} and \eqref{e2} we arrive at \eqref{est_trans}. $\square$
 
\emph{Proof of Theorem 1.} As we deal with a linear system, once we have shown the estimate \eqref{est_trans} 
the rest follows from a standard density argument. Namely, we consider a sequence of smooth functions
$$
H_\epsilon \to H \quad \textrm{in} \quad W^s_p(\Omega).
$$
Let $\sigma_{\epsilon}$ be a solution to \eqref{sys} with $H_{\epsilon}$. 
In particular it satisfies the weak formulation \eqref{weak} with $H_{\epsilon}$. 
As $\sigma_{\epsilon}$
is bounded in $W^s_p$, it converges weakly in $W^s_p$ up to a subsequence (still denoted by $\sigma_{\epsilon}$)
to some $\sigma \in W^s_p$. We have to show that all the integrals in the weak formulation convergence 
but this is obvious due to weak convergence as each integral can be understood as a functional on $W^s_p$. 
We conclude that $\sigma$ is a solution to \eqref{sys}. $\square$     

\smallskip
\noindent
{\bf Acknowledgements.} The work was supported by polish NCN grant UMO-2014/14/M/ST1/00108.

\end{document}